\documentclass{article}
\usepackage{amsmath2000}
\usepackage{amssymb}
\usepackage[dvips]{epsfig}
\usepackage[arrow,matrix]{xy}
\newcommand{\ind}{\operatorname{ind}}
\newcommand{\interior}{\operatorname{int}}
\newcommand{\indt}{{\ind}_\text{\textup{t}}}
\newcommand{\inda}{{\ind}_\text{\textup{a}}}
\newcommand{\Ca}{$C^*$-algebra}

\newcommand{\Oms}{\Omega^{\mathrm{Spin}}}

\newcommand{\fp}{\mathfrak p}
\newcommand{\bR}{\mathbb R}

\newcommand{\bC}{\mathbb C}

\newcommand{\cH}{\mathcal H}
\newcommand{\cC}{\mathcal C}

\newcommand{\bZ}{\mathbb Z}
\newcommand{\co}{\colon\,}
\newcommand{\lp}{\textup{(}}
\newcommand{\rp}{\textup{)}}
\newcommand{\pt}{\text{\textup{pt}}}

\newcommand{\psc}{positive scalar curvature}
\newtheorem{thm}{Theorem}[section]
\newtheorem{cor}[thm]{Corollary}

\newtheorem{prop}[thm]{Proposition}
\newtheorem{remk}[thm]{Remark}
\newtheorem{defin}[thm]{Definition}
\newtheorem{question}[thm]{Question}

\newenvironment{demo}{\par\noindent{\em
Proof}. }{$\square$\par\indent}

\newenvironment{defn}{\begin{defin}\em}{\end{defin}}
\newenvironment{quest}{\begin{question}\em}{\end{question}}
\newenvironment{rem}{\begin{remk}\em}{\end{remk}}
\numberwithin{equation}{section}
\numberwithin{figure}{section}
\begin{document}
\title{Groupoid {\Ca}s and\\ index theory on\\
manifolds with singularities}
\author{Jonathan Rosenberg}
\footnotetext{Mathematics Subject Classification (MSC2000): Primary
58J22. Secondary 19K56, 46L87, 46L85, 46L80.}
\maketitle
\markright{Groupoid {\Ca}s and index theory with singularities}
\begin{abstract}
The simplest case of a manifold with singularities is a manifold $M$
with boundary, together with an identification $\partial M\cong
\beta M\times P$, where $P$ is a fixed manifold. The associated
singular space is obtained by collapsing $P$ to a point. When
$P=\bZ/k$ or $S^1$, we show how to attach to such a space a
noncommutative $C^*$-algebra that captures the extra structure. We
then use this $C^*$-algebra to give a new proof of the Freed-Melrose
$\bZ/k$-index theorem and a proof of an index theorem for manifolds
with $S^1$ singularities. Our proofs apply to the real as well as to
the complex case.  Applications are given to the study of metrics of
\psc. 
\end{abstract}

\section{Introduction}\label{sec:intro}
The simplest example of a ``manifold with singularities'' in the sense
of Sullivan (\cite{Sull},
cf.{\ }also \cite{MorgSull}) and Baas \cite{Baas} is a 
$\bZ/k$-manifold, a manifold with boundary whose boundary
consists of $k$ identical components, all identified with one another.
These were originally introduced for the purpose of giving a geometric
meaning to bordism with $\bZ/k$ coefficients, or to index invariants
with values in  $\bZ/k$, such as the ``signature mod $k$'' (see
\cite{MorgSull}). Later, Freed \cite{FreedZkind} 
and Freed-Melrose \cite{FreedMelrose} were able to give
an analytic version of a index theorem for such manifolds,
and other (or, as some might claim, better) proofs were given by
Higson \cite{HigsonZkind},
Kaminker-Wojciechowski \cite{KamWoj}, and Zhang \cite{Zhang}.
The innovation in Higson's proof was the use of noncommutative
{\Ca}s to model the operator-theoretic part of the index calculation.

However, the ``philosophy of noncommutative geometry'' would suggest
still another approach. Namely, a $\bZ/k$ structure on a manifold
should be modeled by a noncommutative {\Ca}, and then one should
work with this {\Ca} just as one would work with the usual algebra of
functions on a manifold in proving the index theorem. The purpose of
this paper is to implement such an approach, using the
concept of a groupoid {\Ca} as introduced by Renault
(\cite{Renault}, or one can also find a nice exposition in
\cite{Paterson}). This has several
ancillary benefits.  First of all, it treats the $\bZ/k$ index 
just as one treats other kinds of indices in Kasparov-style
index theory, following the outline that one can find in
\cite{Rosrev}. Secondly, it suggests a program for extending everything
to manifolds with other kinds of singularities, of which a 
good general exposition can be found in \cite{Botv}. This
has potential geometric applications as explained for example
in \cite{BotvPSC}.

This paper arose out of lectures I prepared for the Summer
Research Conference on Noncommutative Geometry at Mount  Holyoke
College in June, 2000. It's a pleasure to thank the organizers of that
meeting, Alain Connes, Nigel Higson, and John Roe, as well as the
sponsor of the lectures, the Clay Mathematics Institute, for their
help and support. Some of the results were also presented at the
Workshop on Noncommutative Geometry and Quantization at MSRI in April,
2001.

\section{The {\Ca} of a $\bZ/k$-Manifold}\label{sec:Ca}

First we recall the notion of a (smooth closed) manifold with 
singularities. For
more details, see \cite{Baas}, \cite{Botv}, or \cite{BotvPSC}.
\begin{defn}\label{def:manwithsing}
Let $\Sigma=(P_1,\cdots,P_r)$, where $P_1,\cdots,P_r$
are closed (smooth) manifolds (not necessarily connected). Then a (closed)
\emph{$\Sigma$-manifold}, or \emph{manifold with singularities
$\Sigma$}, means a compact manifold $M$ with boundary, together with a
decomposition $\partial M=\partial_1 M\cup \cdots \cup\partial_r M$.
We require that for each $I=\{i_1,\cdots,i_q\}\subseteq \{1,\cdots,
r\}$, $\partial_I M=\partial_{i_1} M\cap \cdots \cap\partial_{i_q} M$
is a manifold with boundary
\[
\partial\left(\partial_I M\right)=\bigcup_{j\notin I}
\partial_{I\cup\{j\}} M
\]
equipped with a  product identification
\[
\phi_I\co \partial_I M \stackrel{\cong}{\longrightarrow} 
\beta_I M \times P^I,
\]
where $\beta_I M$ is a manifold and
$P^I=P_{i_1} \times \cdots \times P_{i_q}$. In addition, we
require compatibility of these maps; i.e., if $J\subset I$
and $\iota_{IJ}\co \partial_I M \to \partial_J M$ denotes the
inclusion, then 
\[
\xymatrix{
\beta_I M \times P^I 
\cong \beta_I M \times P^J\times P^{I\smallsetminus J}
\ar[rr]^(.65){\phi_J\circ\iota_{IJ}\circ\phi_I^{-1}}
\ar[rd]_{\text{proj}_2}
& & \beta_J M \times P^J
\ar[ld]^{\text{proj}_2}\\
&P^J&
}
\]
must commute.
\end{defn}

We will mostly be concerned with the easy case where $r=1$
and $P_1=P$, a single fixed closed manifold. Then the
conditions amount to saying that $M$ is a compact manifold with an
identification $\partial M \stackrel{\cong}{\longrightarrow} 
\beta M \times P$. Perhaps the most interesting special
case is still the original situation of a $\bZ/k$-manifold,
where $P$ is the disjoint union of $k$ points. This 
case is illustrated in Figure \ref{fig:Z3man}, with $k=3$.

\setcounter{figure}{0}
\begin{figure}[htb]\refstepcounter{figure}\label{fig:Z3man}
\begin{center}\epsfig{file=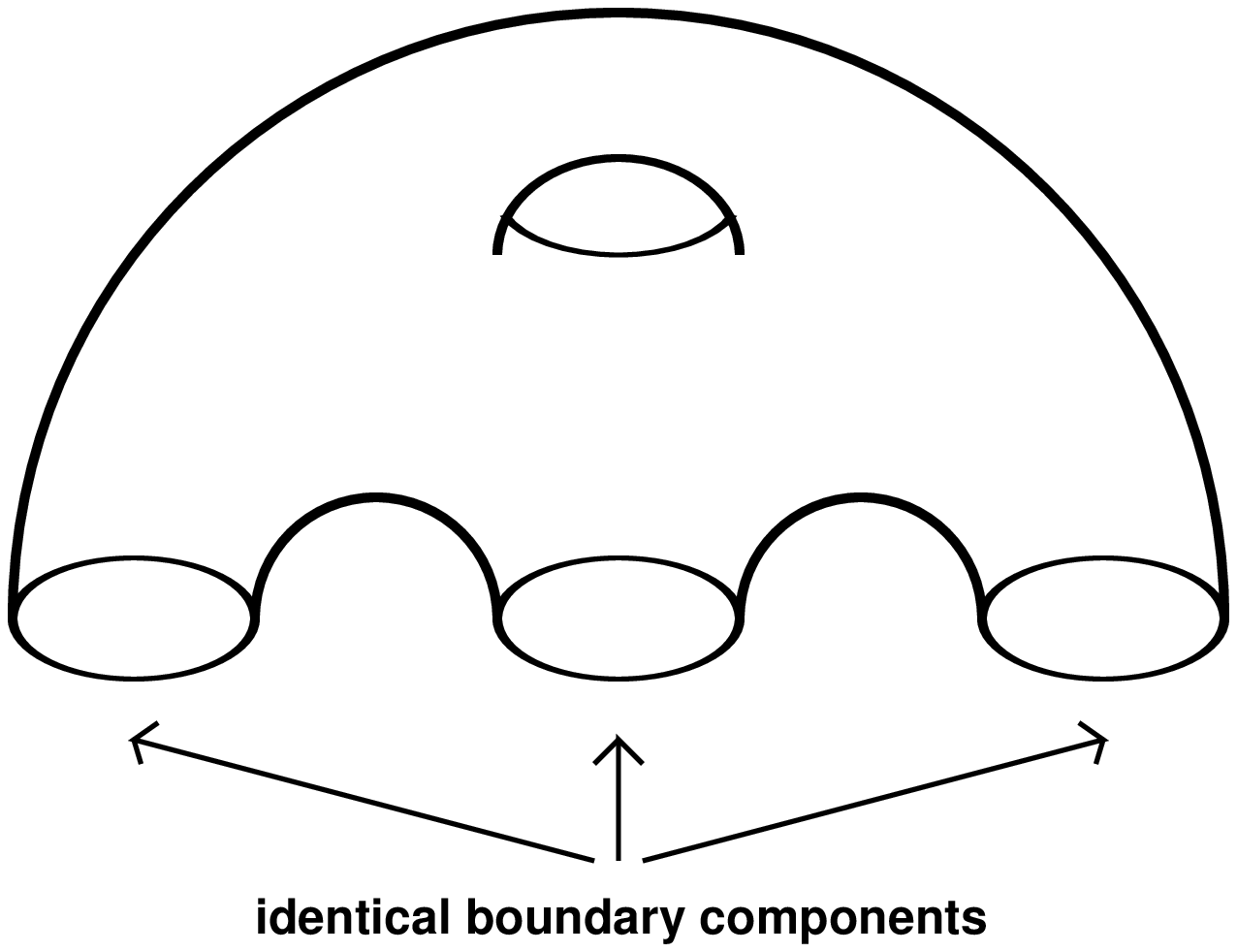,%
width=3.5in}\setcounter{figure}{0}
\caption{A $\bZ/3$-manifold}\end{center}\end{figure}

Note that while $M$ itself is smooth, there is a singular
space $M_\Sigma$ associated to it, obtained by collapsing
each $\partial_I M$ to $\beta_I M$ via the composition
\[
\xymatrix{
M_I \ar[r]^(.4){\phi_I} & \beta_I M \times P^I \ar[d]^{\text{proj}_1}\\
& \beta_I M.
}
\]
The compatibility condition above guarantees that the resulting
quotient space $M_\Sigma$ is well-defined. In the easy case where $r=1$
and $P_1=P$, $M_\Sigma$ is $M$ with $\partial M 
\stackrel{\cong}{\longrightarrow} \beta M \times P$ collapsed
to $\beta M$ by collapsing each $\{x\}\times P$, $x\in \beta M$, 
to a point. In all cases, the map $M\twoheadrightarrow M_\Sigma$ is
injective on the interior of $M$.

Roughly speaking, the \Ca, $C^*(M;\Sigma)$, of a manifold 
$M$ with singularities $\Sigma$ should be the algebra of functions 
on the singular quotient space $M_\Sigma$. However, this is too 
coarse an invariant, as it
doesn't take all the extra structure into account. Instead,
since $M_\Sigma$ is the quotient of $M$ by an equivalence relation
$\sim$ encoding the $\Sigma$-structure on $M$, we want to
use something like
the {\Ca} of the equivalence relation. Here we are following
the philosophy of noncommutative geometry as expounded in
\cite{Connes}: points in the same equivalence class should
``talk to each other'' but not be collapsed to one.  However, this is 
still not exactly right, for a number of reasons:
\begin{enumerate}
\item Later we will want to do analysis without having to
worry about the boundary. For this reason (following a trick
in \cite{HigsonZkind}) we add half-infinite
cylinders onto the boundary first, as in Figure 
\ref{fig:Z3manextended}.
\item We need to construct the {\Ca} so that it has the ``right''
$K$-theory. For instance, in the case of a $\bZ/k$-manifold,
$K_*(C^*(M;\Sigma))$ should be related to $K$-theory of $M$ with
$\bZ/k$ coefficients.
\item We want to construct $C^*(M;\Sigma)$ as the {\Ca} of
a suitable locally compact groupoid, in the sense of \cite{Renault}.
Such a groupoid is not so hard to construct in the case of a 
$\bZ/k$-manifold (see Definition \ref{def:CaZk}), but in the case of
general singularities, it is not clear how to proceed.
\end{enumerate}

\begin{figure}[htb]\refstepcounter{figure}\label{fig:Z3manextended}
\begin{center}\epsfig{file=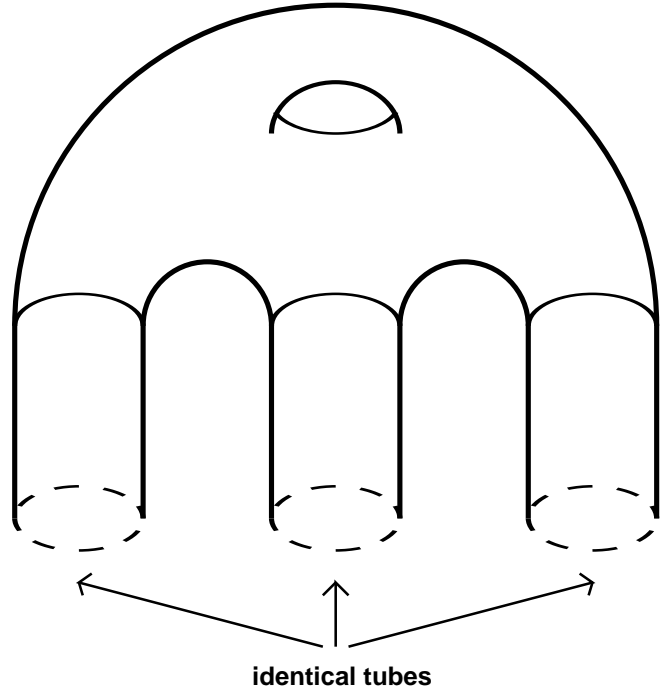,%
width=3.5in}
\setcounter{figure}{1}
\caption{A $\bZ/3$-manifold with cylinders added}
\end{center}\end{figure}

\begin{defn}\label{def:CaZk}
Let $k\ge 2$ and let $M$ be a compact $\bZ/k$-manifold, that is a 
manifold with singularities $\Sigma=(\bZ/k)$ in the sense of
Definition \ref{def:manwithsing}. Recall that this means that,
as part of the given structure on $M$, we have a diffeomorphism
$\phi\co \partial M\to \beta M\times \bZ/k$, for
some closed manifold $\beta M$. Let $N=M\cup_{\partial M}
{\partial M} \times [0,\infty)$, where $\partial M\subset M$
is identified to $\partial M\times \{0\}$. This is a manifold without
boundary, usually non-compact. Note that $N$ is homeomorphic to
the interior of $M$, i.e., to $M\smallsetminus\partial M$,
via the collaring theorem. We define an equivalence relation
$\sim$ on $N$ as follows. Points in $M$ (including points in $\partial M$)
are equivalent only to themselves. A point $(x,t)\in {\partial M} 
\times (0,\infty)$ is equivalent to a point $(y,s)\in {\partial M} 
\times (0,\infty)$ if and only of $t=s$ and $p(x)=p(y)$, where
$p\co \partial M\to \beta M$ is $\phi\co \partial M\to \beta M\times
\bZ/k$ followed by projection onto the first factor.

Let $G\subset N\times N$ be the equivalence relation $\sim$ viewed as
a groupoid. We observe that $G$ is locally closed in $N\times N$.
Indeed, the closure $\overline G$ of $G$ in $N\times N$ is easily seen
to consist of $G\cup G'$, where $G'$ is the equivalence relation
on $\partial M$ identifying points which project to the same point
in $\beta M$. So 
\begin{multline*}
\overline G\smallsetminus G\cong \left\{(x,i,j)\mid
x\in \beta M,\ 1\le i,\,j\le k,\ i\ne j \right\}\\
\cong \beta M \times ((k^2-k)\text{ points})
\end{multline*}
is compact and thus $G$ is open in $\overline G$, so that $G$
is locally compact in the relative topology from $N\times N$.
(Note that this argument would break down in the case where 
$\Sigma=(P)$, $\dim P>0$.) The unit space of the groupoid $G$ is of
course $N$, and the range map $r\co G\to N$ is a local homeomorphism,
since this is obvious over the interior of $M$ and over $\partial M
\times (0,\infty)$, whereas the only points in $G$ over 
$\partial M\times\{0\}$ are of the form $(x,x)$, $x\in \partial M$,
whose small neighborhoods in $G$ have projection under $r$ that
miss all but one of the cylinders $\beta M \times\{j\} \times
[0,\infty)$. Hence by \cite{Renault}, Proposition 2.8, $G$
has a Haar system which is essentially unique, and the {\Ca}
$C^*(G)$ is well-defined. We denote it by $C^*(M;\bZ/k)$.
Note, incidentally, that in the construction of $C^*(M;\bZ/k)$,
one can use either real or complex scalars. When it is necessary
to distinguish the real and complex {\Ca}s, we will denote the
former by $C^*_\bR(M;\bZ/k)$ for emphasis.
\end{defn}
\begin{prop}\label{prop:Zkalgstructure}
Let $k\ge 2$ and let $M$ be a compact $\bZ/k$-manifold.
Then $C^*(M;\bZ/k)$ fits into a short exact sequence of {\Ca}s
\begin{equation}\label{eq:Zkalgstructure}
0\to C_0(\bR)\otimes C(\beta M)\otimes M_k \to C^*(M;\bZ/k)\to C(M) \to 0.
\end{equation}
This is valid with either real or complex scalars.
\end{prop}
\begin{demo}Note that $M$ is closed in $N$ and on $N$ the equivalence
relation $\sim$ of Definition \ref{def:CaZk} is trivial. So we obtain
a quotient {\Ca} of $C^*(M;\bZ/k)$ isomorphic to $C(M)$, and the
kernel of the quotient map $C^*(M;\bZ/k)\to C(M)$ must be attached to
the complementary open set, $\partial M\times (0,\infty)$. Since the
inverse image of this set in $G$ splits as a product $\beta M\times
(0,\infty) \times (\bZ/k)\times (\bZ/k)$, and the {\Ca} of
the groupoid
$(\bZ/k)\times (\bZ/k)$ is just $M_k$, the result follows.
\end{demo}
\begin{prop}\label{prop:ZkalgKthy}
Let $k\ge 2$ and let $M$ be a connected compact $\bZ/k$-manifold.
Then the connecting map in the long exact $K$-theory sequence
associated to the extension \textup{(\ref{eq:Zkalgstructure})}
may be identified with the map $\iota^*\co K^*(M)\to
K^*(\beta M)$, followed by multiplication by $k$. Here $\iota\co \beta M \to
M$ denotes the inclusion.  If complex scalars are used, $K^*$ means 
$KU^*$, and if real scalars are used, $K^*$ means $KO^*$.
\end{prop}
\begin{demo}The connecting map of \textup{(\ref{eq:Zkalgstructure})}
in $K$-theory is
\begin{multline}\label{eq:ZkalgKthy}
K_i(C(M))\cong K^{-i}(M) \to \\
K_{i-1}\bigl(C_0(\bR)\otimes C(\beta M)\otimes
M_k\bigr) \cong K_{i-1}\bigl(C_0(\bR)\otimes C(\beta M)\bigr)\\
 \cong K^{-i+1}(\bR\times \beta M)\cong K^{-i}(\beta M).
\end{multline}
Now let $C^*(\partial M;\bZ/k)$ be defined like $C^*(M;\bZ/k)$,
in the sense that we use an equivalence relation on $\partial M
\times [0,\infty)$ which is trivial on $\partial M\times \{0\}$
and identifies points in $\partial M
\times (0,\infty)$ having the same image in $\beta M\times 
(0,\infty)$. Then in analogy with extension (\ref{eq:Zkalgstructure})
we have an extension
\begin{equation}\label{eq:Zksmallalgstructure}
0\to C_0(\bR)\otimes C(\beta M)\otimes M_k \to C^*(\partial M;\bZ/k)
\to C(\partial M) \to 0
\end{equation}
in which $C(\beta M)$ splits off as a tensor factor.  In other
words, $C^*(\partial M$; $\bZ/k)\cong C(\beta M)\otimes
C^*(\bZ/k;\bZ/k)$, where we have an extension
\begin{equation}\label{eq:ZkalgZkstructure}
0\to C_0(\bR)\otimes M_k \to C^*(\bZ/k;\bZ/k) \to 
C(\pt)^k\to  0.
\end{equation}
Now because of the commutative diagram
\[
\xymatrix@C-1.8pt{
0\ar[r]\ar@{=}[d] & C_0(\bR)\otimes C(\beta M)\otimes M_k 
\ar[r]\ar@{=}[d] & C^*(M;\bZ/k) \ar[r] \ar@{{>>}}[d] 
& C(M) \ar[r] \ar@{{>>}}[d] & 0 \ar@{=}[d] \\
0\ar[r] & C_0(\bR)\otimes C(\beta M)\otimes M_k \ar[r] &
C^*(\partial M;\bZ/k)\ar[r] & C(\partial M) \ar[r] & 0 \\
}
\]
the map of (\ref{eq:ZkalgKthy}) factors through the restriction map
$K^{-i}(M)\to K^{-i}(\partial M)$, followed by
connecting map for (\ref{eq:Zksmallalgstructure}),
which in turn is the external product of the identity map
on $ K^{-i}(\beta M)$ with the connecting map for
(\ref{eq:ZkalgZkstructure}).  But $K^*(\partial M)\cong
K^*(\beta M)^k$, and the map $K^*(M)\to K^*(\partial M)$
is just $\iota^*\co K^*(M)\to K^*(\beta M)$, followed by 
the diagonal inclusion of $K^*(\beta M)$ into a product of $k$
copies of itself. So we can rewrite the map of (\ref{eq:ZkalgKthy}) 
as the composite of $\iota^*\co K^*(M)\to K^*(\beta M)$ with the
external product of the identity map
on $ K^*(\beta M)$ with the connecting map for the result of
collapsing the $k$ copies of the scalars in 
(\ref{eq:ZkalgZkstructure}) to one:
\begin{equation}\label{eq:Zkalgptstructure}
0\to C_0(\bR)\otimes M_k \to C^*(\pt;\bZ/k) \to 
C(\pt)\to  0.
\end{equation}

We also only need to compute the connecting map in $K_0$. (That's
because the  connecting map for (\ref{eq:Zkalgptstructure}) is easily
seen to be a map of $K_*(\pt)$-modules from $K_*(\pt)$
to itself, and so it's determined by what happens to the generator
in degree $0$.) Now one computes that
\[
C^*(\pt;\bZ/k)=
\{f\in C_0([0,\infty), M_k) \mid f(0) \text{ a multiple of } I_k\}.
\]
This is simply the mapping cone of the inclusion of
the scalars into $M_k$ as multiples of the $k\times k$ identity
matrix, so the connecting map is multiplication by $k$, as
required.
\end{demo}
\begin{rem}In fact Propositions \ref{prop:Zkalgstructure} and
\ref{prop:ZkalgKthy} are also valid when $k=1$, but in this case,
$\partial M=\beta M$ and $C^*(M;\bZ/k)$ is simply $C_0(N)$.
Since $N$ is properly homotopy equivalent to the interior of $M$, its
$K$-theory is just the relative $K$-theory of the pair $(M,\partial
M)$, and so all statements are obvious in this case.
\end{rem}
\begin{cor}\label{cor:KthyZkpt} 
The $K$-theory of the {\Ca} $C^*(\mathrm{pt};\bZ/k)$
defined in \textup{(\ref{eq:Zkalgptstructure})} is just
$K_i(C^*(\mathrm{pt};\bZ/k))\cong K^{-i-1}(\mathrm{pt};\bZ/k)$. 
The dual theory
{\lp}often known as  ``K-homology,'' though it is a \emph{cohomology}
theory on {\Ca}s{\rp} is given by $K^{-i}(C^*(\mathrm{pt};\bZ/k))
\cong K_i(\mathrm{pt};\bZ/k)$.
\end{cor}
\begin{demo}This is obvious from Proposition \ref{prop:ZkalgKthy}, the
long exact sequences, and the universal coefficient theorem.
\end{demo}

\section{The {\Ca} of an $\eta$-Manifold}\label{sec:eta}

As mentioned above, we are not sure how to give a good definition of
$C^*(M; \Sigma)$ in the case of general singularities. However, there
is one case which is tractable and interesting, the case of an
$\eta$-manifold. 
\begin{defn}\label{def:eta}
An \emph{$\eta$-manifold} (cf.{\ }\cite{BotvPSC})
is a manifold with singularities
$\Sigma=(\eta)$ in the sense of Definition \ref{def:manwithsing}, where
$\eta$ denotes $S^1$ with its non-bounding spin structure.  If we
forget the spin structure, an $\eta$-manifold is just a manifold $M$ with
singularities $\Sigma=(S^1)$, but later we will require $M$  to have a
spin structure inducing a product spin structure $(\beta M, s)\times
\eta$ on its boundary. (Here $s$ is a spin structure on $\beta M$.) 
As is well known, $\eta$ generates $\Oms_1 \cong \bZ/2$, and is the
image of the generator of $\pi_1^s\cong \bZ/2$ with the same name.
\end{defn}

\begin{defn}\label{def:Caeta}
Let $M^n$ be an $\eta$-manifold in the sense of Definition
\ref{def:eta}. Thus $M$ is a compact manifold with boundary, together
with a diffeomorphism $\phi\co \partial M\to \beta M\times S^1$
(and some spin structure data). Let $N=M\cup_{\partial M}
{\partial M} \times [0,\infty)$, where $\partial M\subset M$
is identified to $\partial M\times \{0\}$. This is a manifold without
boundary, usually non-compact. Let $N'$ be the quotient space of $N$
(locally compact, but not a manifold) obtained by collapsing
$\partial M\subset M$ to $\beta M$ via the projection map
$\beta M\times S^1 \to \beta M$. Note that the algebra $C_0^\bR(N')$
can be identified with the subalgebra of functions $f\in C_0^\bR(N)$ which
are constant along $\{x\}\times S^1$ for each $x\in \beta M$.
Identify $S^1$ with the unit circle in the complex plane and let 
the group $S^1$ act on $N'$ as follows: the action is trivial on 
the interior of $M$ and on the image of $\partial M$, and $S^1$ acts on
\[
{\partial M} \times (0,\infty) \cong \beta M\times S^1 \times
[0,\infty) 
\]
by the formula $z\cdot(x,w,t)=(x,zw,t)$, $x\in\beta M$,
$z,w\in S^1$, $t\in\bR$. Alternatively, the image of 
${\partial M} \times [0,\infty)$ in $N'$ can be identified with
$\beta M\times \bC$, with a product action
where $S^1$ acts trivially on $\beta M$ and
on $\bC$ by multiplication.  (Here $\beta M\times \{0\}$ corresponds
to the image of ${\partial M} \times  \{0\}$, and 
$\beta M\times \{z\in\bC \mid \vert z\vert =t\}$
is identified with ${\partial M} \times \{t\}$ when $t>0$.)
The fixed point set for the $S^1$-action on $N'$ is the singular
space $M_\Sigma$ ($M$ with $\partial M\cong
\beta M\times S^1$ collapsed to $\beta M$).

Let $C^*_\bR(M;\eta)$ and $C^*(M;\eta)$ be the real and complex
transformation {\Ca}s $C^\bR_0(N')\rtimes S^1$ and
$C_0(N')\rtimes S^1$. Also let $C^*_\bR(\pt;\eta)$ denote
$C_0^\bR(\bC)\rtimes S^1$, and similarly with $C^*(\pt;\eta)$. 
We have a canonical $S^1$-equivariant
proper ``collapse'' map $c\co N'\to \bC$ sending $M$ to $0$ and
projecting 
\[
\beta M\times \{z\in\bC \mid \vert z\vert =t\} \to \{z\in\bC \mid
\vert z\vert =t\} 
\]
for $t>0$, so this induces a dual map of {\Ca}s
$C^*_\bR(\pt;\eta) \to C^*_\bR(M;\eta)$, and similarly in the
complex case.
\end{defn}
\begin{prop}\label{prop:etaalgKthy}The $S^1$-equivariant real
$K$-theory of $\bC$ with compact supports {\lp}for the action of $S^1$ by 
multiplication{\rp} may be identified with the $K$-theory of
$C^*_\bR(\text{\textup{pt}};\eta)$, and is given by
\[
KO_i(C^*_\bR(\text{\textup{pt}};\eta))\cong KO^{-i}_{S^1}(\bC)
\cong KO^{-i-2}_{S^1}(\text{\textup{pt}})
\]
for all $i$ {\lp}though one has periodicity mod $8${\rp}. 
These groups are computed in \textup{\cite{AS}}. The {\lp}Kasparov{\rp}
dual theory is given by
\[
KO_{i}^{S^1}(\bC)\cong KO_{i-2}^{S^1}(\text{\textup{pt}})
\]
for all $i$ {\lp}again with periodicity mod $8${\rp}. Similarly,
$K_{i}^{S^1}(\bC)\cong R(S^1)=\bZ[t,\,t^{-1}]$ for $i$ even, $0$
for $i$ odd.
\end{prop}
\begin{demo}The identification of equivariant $K$-theory (for a
compact group action) with $K$-theory of the crossed product may be
found in \cite{Julg}. The proof is only stated for the complex case,
but everything goes over the real case as well.
The rest of this is equivariant Bott periodicity, for which
one can see \cite{Segal} for the case of $K$-cohomology,
and \cite{Kas}, \S5, for the case of Kasparov $K$-homology.
\end{demo}

\section{Three index theorems}\label{sec:index}

\begin{defn}\label{def:singop}
Let $M$ be a closed $\Sigma$-manifold, where $\Sigma=(P)$, $P=\bZ/k
\text{ or }\eta$. Recall that $M$ comes with a diffeomorphism $\phi\co
\partial M\to \beta M\times P$. For purposes of this section, an
elliptic operator $D$ on $M$ will mean an elliptic pseudodifferential
operator which on a small collar neighborhood of the boundary, 
$\partial M\times [0,\varepsilon)$, is invariant under translation
normal to the boundary (i.e., is the restriction of a $\bR$-invariant
elliptic operator on $\partial M\times \bR$) and is also invariant
under translation in $P$. (Recall that in our case, $P$ is a
compact group, either $\bZ/k$ or $S^1$, and we can identify $\partial
M\times [0,\varepsilon)$ with $\beta M\times P\times [0,\varepsilon)$
via $\phi$.)  The (graded) vector bundle on which such
an operator acts must satisfy a similar invariance condition.
The primary examples of such operators are the
\emph{standard elliptic operators}: the Euler characteristic operator
of a Riemannian manifold, the signature operator of an oriented
Riemannian manifold, the Dirac operator of a spin (or
spin$^c$) Riemannian manifold, 
and the Dolbeault operator of a K\"ahler manifold.  In all cases,
the associated Riemannian metric has to be a product metric
(coming from a metric on $\beta M$ and from the standard metrics on $P$
and $\bR$) in a neighborhood of the boundary.
\end{defn}
\begin{defn}\textbf{(the analytic $\bZ/k$-index)}\label{def:analZkind}
Let $(M, \phi\co\partial M\stackrel{\cong}{\to}
\beta M\times \bZ/k)$ be a closed
$\bZ/k$-manifold, and let $D$ be an elliptic operator on $M$ in the
sense of Definition \ref{def:singop}. We assume $D$ is acting on
sections of some vector bundle, equipped with a grading and suitable
Clifford module data so that $D$ locally defines a class in $K_i$
for some $i$, where $K_*$ denotes either real or complex $K$-homology, as
appropriate. The main examples we have in mind are the following:
\begin{itemize}
\item $M^i$ is a spin manifold and $(\beta M)^{i-1}$ has the induced spin
structure. $D$ is the $\cC\ell_i$-linear real Dirac operator in the sense
of \cite{LawM}, Chapter II, \S7, and
Chapter III, \S16. ($\cC\ell_i$ denotes the Clifford algebra of $\bR^i$.)
In this case $D$ locally defines a class in $KO_i$. 
\item $M^{2n}$ is a spin$^c$ manifold and $(\beta M)^{2n-1}$ has the 
induced spin$^c$ structure. $D$ is the  complex Dirac operator acting
on the spinor bundle with $\bZ/2$-grading defined by the half-spinor
bundles. In this case $D$ locally defines a class in $K_0$. 
\item $M^{2n}$ is an oriented even-dimensional manifold and $D$ is the
signature operator as in \cite{LawM}, Chapter II, \S6, which locally
defines a class in $K_0$.
\item $M^{2n}$ is a complex manifold (of complex dimension $n$) with a
K\"ahler metric which near the boundary is compatible with the
$\bZ/k$-structure, and $D$ is the Dolbeault operator. In this case $D$ 
locally defines a class in $K_0$. 
\end{itemize}
Since, in a neighborhood of the boundary, $D$ is invariant under
translations normal to the boundary of $M$, it naturally extends to an
elliptic operator on the manifold $N$ of Definition \ref{def:CaZk}.
We claim that $D$ defines a class $[D]\in
K^{-i}(C^*(M;\bZ/k))=KK(C^*(M;\bZ/k), \cC\ell_i)$. Indeed, this is
obvious from the fact that $D$ defines a $C_0(N)$-$\cC\ell_i$ Kasparov
bimodule which is equivariant for the groupoid $G$ (in Definition
\ref{def:CaZk}), and thus gives a Kasparov bimodule for the groupoid
$C^*$-algebra. Since a homotopy of metrics gives rise to a homotopy
of operators, the Kasparov class $[D]$ is independent of the choice
of metric. Note that via the open inclusions $\bR\times\beta M
\hookrightarrow C^*(M;\bZ/k)^{\widehat{}}$ and $\text{int}\,
M\hookrightarrow C^*(M;\bZ/k)^{\widehat{}}$, $[D]$ restricts to the
usual class defined by the relevant elliptic operator on the open
manifold $\bR\times\beta M$ or $\text{int}\,M$.
The \emph{analytic $\bZ/k$-index} $\inda(D)$ of $D$ is defined to be the
image of the class $[D]$ under the map $c_*\co K^{-i}(C^*(M;\bZ/k))
\to K^{-i}(C^*(\pt;\bZ/k))\cong K_i(\pt;\bZ/k)$. The
map $c_*$ is dual to the inclusion 
\begin{equation}
\label{eq:collapse}
c^*\co C^*(\pt;\bZ/k)\hookrightarrow
C^*(M;\bZ/k), 
\end{equation}
or equivalently can be viewed as the map in
``$K$-homology'' induced by the ``collapse map'' $c$ of $M$ to a
point. The identification of $K^{-i}(C^*(\pt$; $\bZ/k))$ comes from
Corollary \ref{cor:KthyZkpt}.
\end{defn}
\begin{defn}\textbf{(the topological
$\bZ/k$-index)}\label{def:topZkind}
Let notation be as in Definition \ref{def:analZkind} above. 
Let $[\sigma(D)]\in K^*(T^*\!M)$ or $KR^*(T^*\!M)$, the $K$-theory
with compact supports of the cotangent bundle of $M$, be the class of
the principal symbol of the operator. In the real case we need to view
$T^*\!M$ as a Real space, as explained in \cite{AtKR}, or in 
\cite{LawM}, Chapter III, \S16. Note that $[\sigma(D)]$
is invariant under the identifications on the boundary, i.e., it comes
by pullback from the quotient space $T^*\!M_\Sigma$ (the image of $T^*\!M$
with the $k$ copies of $T^*\!M|_{\beta M}$ collapsed to one) under the
collapse map $M\twoheadrightarrow M_\Sigma$. Following
\cite{FreedZkind} we define the \emph{topological $\bZ/k$-index}
$\indt D$ of
$D$ as follows.  Start by choosing an embedding $\iota\co(M,\partial M)
\hookrightarrow (D^{2r},
S^{2r-1})$ of $M$ into a ball of sufficiently large even dimension $2r$
($2r$ divisible by $8$ in the real case), for which $\partial
M$ embeds $\bZ/k$-equivariantly into the boundary (if we identify
$S^{2r-1}$ with the unit sphere in $\bC^r$, $\bZ/k$ acting as usual by
multiplication by roots of unity). In the complex case, with $D$
anti-commuting with a $\bZ/2$-grading so as to give a class in $K_0$,
we take the push-forward map on complex $K$-theory
\[
\iota_!\co 
K^0(T^*\!M)\to \widetilde{K}^0(T^*\!D^{2r})
\cong \widetilde{K}^0(D^{2r})
\]
and observe that $\iota_!([\sigma(D)])$ descends to 
$\widetilde{K}^0(M_k^{2r})\cong K^0(\pt;\bZ/k)\cong \bZ/k$,
$M_k^{2r}$ the Moore space obtained by
dividing out by the $\bZ/k$-action on the boundary of $D^{2r}$,
and call the image the topological index of $D$,
$\indt(D)$. The real case is
similar, except that we need to use the push-forward map on Real
$K$-theory instead, getting a topological index in
$KO^*(\pt;\bZ/k)$.
\end{defn}
\begin{thm}[cf{.\ }\cite{FreedZkind}, \cite{FreedMelrose},
\cite{HigsonZkind}, and \cite{KamWoj}]
\label{thm:Zkindthm}
Let $(M, \phi\co\partial M\stackrel{\cong}{\to}
\beta M\times \bZ/k)$ be a closed
$\bZ/k$-manifold, and let $D$ be an elliptic operator on $M$ in the
sense of Definition \ref{def:singop}. Then the analytic index of
$D$ in $K_i(\text{\textup{pt}};\bZ/k)$, in the sense of Definition
\ref{def:analZkind} coincides with the topological index $\indt D$
of $D$ in the
sense of Definition \ref{def:topZkind}. {\lp}This is valid
in both the real and complex cases.{\rp}
\end{thm}
\begin{demo}First consider the complex case. In the only interesting
case, $M$ is even-dimension and $D$ anti-commutes with a
$\bZ/2$-grading, so as to locally define a class in $K^0$. 
Let $N\cong \text{int}\,M$ be as in Definition \ref{def:CaZk}. By the
Kasparov-theoretic proof of the usual index theorem (\cite{Black},
Chapter IX, \S24.5), the class of $D$ in $K_0(N)\cong K_0(\text{int}\,
M)$ is the Kasparov product $[\sigma(D)]\widehat\otimes_{C_0(T^*\!M)}
\alpha$, where 
\[
\alpha\in KK\bigl(C_0(T^*\!M)\otimes C_0(\text{int}\, M),\, 
\bC\bigr)
\]
is the canonical class coming from the Dolbeault complex 
$\overline\partial$ of the
canonical almost complex structure on $T^*\!M$, the Thom isomorphism, 
and the projection map $T^*\!M\to M$. (In Blackadar's book
the proof is given in the case where $\partial M=\emptyset$, but the
case where $M$ has a boundary works the same way, once one notices
that $D$ and $\overline\partial$
define Kasparov bimodules for $C_0(\text{int}\, M)$,
though not for $C(M)$, since we have not imposed any boundary
conditions.) Now one can observe that everything in sight is
compatible with the $\bZ/k$-structure, and so descends to the groupoid
algebra. In other words, with notation as in Definition
\ref{def:topZkind}, we have
\[
[D]=[\sigma(D)]\widehat\otimes_{C_0(T^*\!M_\Sigma)}
\widehat\alpha \in K^0(C^*(M;\bZ/k)),
\]
where now
\[
\widehat\alpha\in KK\bigl(C_0(T^*\!M_\Sigma)\otimes C^*(M;\bZ/k),\, 
\bC\bigr)
\]
is the groupoid-equivariant version of $\alpha$, and we now view
$[\sigma(D)]$ as living in $K^0(T^*\!M_\Sigma)$.

But now $c_*\co K^0(C^*(M;\bZ/k))
\to K^0(C^*(\pt;\bZ/k)) \stackrel{\cong}{\longrightarrow}
K_0(\pt;\bZ/k)$ can be
viewed as Kasparov product with the homomorphism $c^*$ of equation
(\ref{eq:collapse}). 
So by associativity of the Kasparov product, we
compute that
\[
\inda(D)=[c^*]\widehat\otimes_{C^*(M;\bZ/k)} [D]
=[\sigma(D)]\widehat\otimes_{C_0(T^*\!M_\Sigma)}\Bigl(
[c^*]\widehat\otimes_{C^*(M;\bZ/k)}\widehat\alpha \Bigr).
\]
So we just need to identify the right-hand side of this equation with
$\indt(D)$. However, by Definition \ref{def:topZkind}
$\indt(D)={\widehat\iota}_!([\sigma(D)])$, where
\[
{\widehat\iota}_!\co K^0(T^*\!M_\Sigma)\to K^0(T^*\!D^{2r}_\Sigma)
\cong K^0(M^{2r}_k)
\]
is the push-forward map on $K$-theory. And examination of the
definition of ${\widehat\iota}_!$ shows it is precisely the Kasparov
product with 
\[
[c^*]\widehat\otimes_{C^*(M;\bZ/k)}\widehat\alpha,
\]
followed by Kasparov product with a ``Poincar\'e duality'' element
\[
\delta \in KK(C(\pt), C^*(\pt;\bZ/k)\otimes C_0(M_k^{2r}))
\]
implementing the isomorphism $K^0(C^*(\pt;\bZ/k)) 
\stackrel{\cong}{\longrightarrow} K_0(\pt;\bZ/k)$.
The proof in the real case follows exactly the same outline, except
that one has to use the Real structure of the cotangent bundle, i.e.,
replace $C_0(T^*\!M)$ by $\{f\in C_0(T^*\!M)\mid
f(\tau(x))=\overline{f(x)} \}$, where $\tau$ is the involution on
$T^*\!M$ that is multiplication by $-1$ on each fiber.
\end{demo}

\begin{defn}\textbf{(the analytic $\eta$-index)}\label{def:analetaind}
Let $(M, \phi\co\partial M\stackrel{\cong}{\to}
\beta M\times \eta)$ be a closed $\eta$-manifold. We let $i$ be the
dimension of $M$ and assume $M$ is equipped with a spin structure
inducing a product spin structure on $\partial M^{i-1}\cong 
\beta M^{i-2}\times S^1$,
coming from some spin structure on $\beta M$ and the
\emph{non-bounding} spin structure on $S^1$. We give $M$ a Riemannian
structure compatible with its $\Sigma$-structure. 
For simplicity, we'll first suppose
$i=\dim M$ is even  and consider $D_E$, the complex Dirac operator on $M$
with coefficients in some auxiliary bundle $E$ (whose restriction to 
$\partial M$ is pulled back from $\beta M$).
The operator $D_E$ is defined using a Hermitian connection on $E$
compatible with the $\eta$-manifold structure near $\partial E$.
Since $D_E$ is translation-invariant in a
neighborhood of $\partial M$, it extends to an operator on the open
manifold $N=M\cup_{\partial M}\partial M\times [0,\infty) \cong
\interior M$. Thus we have a Kasparov class $[D_E]\in K_0(N)$. (Since
$N$ is non-compact, this is to be interpreted as $KK(C_0(N),\bC)$.)
Let $c\co N\to\bC$ be the ``collapse'' map collapsing $M$ to $0\in\bC$
and sending $\beta M\times S^1 \times [0,\infty)$ first to
$S^1 \times [0,\infty)$ (by collapsing the $\beta M$ factor) and then
to $\bC$ by means of ``polar coordinates'' ($(e^{i\theta},t)\mapsto
te^{i\theta}$).
We define the \emph{analytic $\eta$-index} $\inda(D_E)$ of $D_E$ to be 
$c_*([D_E])$, or in other words the Kasparov product of
$[D_E]\in KK(C_0(N), \bC)$ with the class of the homomorphism
\[
c^*\co C_0(\bC)\hookrightarrow C_0(N).
\]
It takes its values in $K_0(\bC)\cong\bZ$.

Next, we consider the case of the $\cC\ell_i$-linear real Dirac
operator with coefficients in a real vector bundle $E$, whose
restriction to $\partial M$ again comes from a bundle $E_1$ over
$\beta M$. If we use real instead of complex $K$-theory, the same
procedure as in the complex case gives a $KO$-valued analytic index
$\inda(D_E)\in KO_i(\bC)\cong KO_{i-2}(\pt)$.
\end{defn}
\begin{thm}[$\eta$-manifold index theorem]\label{thm:etaindthm}
Let $(M, \phi\co\partial M\stackrel{\cong}{\to}
\beta M\times \eta)$ be a closed spin {\lp}or spin$^c${\rp}
$\eta$-manifold of even dimension, and let $E$
be a vector bundle on $M$ whose restriction to $\partial M$ is
pulled back from a bundle $E_1$ on $\beta M$. Fix a Riemannian
structure on $M$ compatible with its $\Sigma$-structure, and let $D_E$
be the Dirac operator on $M$ with coefficients in $E$, computed with
respect to a Hermitian connection on $E$ whose restriction to a
neighborhood of the boundary is pulled back from
$\beta M$. Then $\inda(D_E)=
\ind\left({\left(D_{\beta M}\right)}_{E_1}\right)$, the index 
{\lp}in complex $K$-theory{\rp} of
the Dirac operator on $\beta M$ with coefficients in $E_1$, which in
turn is computed by applying the usual Atiyah-Singer Theorem on $\beta M$.
\end{thm}
\begin{demo}Consider the following commutative diagram with exact
rows:
\[
\xymatrix{
0\ar[r]\ar@{=}[d] &C(\beta M)\otimes C(S^1)\otimes  C_0(\bR)
\ar[r]\ar@{{>>}}[d] & C_0(N) \ar[r] \ar@{{>>}}[d] 
& C(M) \ar[r] \ar@{{>>}}[d] & 0 \ar@{=}[d] \\
0\ar[r] & C_0(\bC\smallsetminus \{0\}) \ar[r] &
C_0(\bC)\ar[r] & C(\{0\}) \ar[r] & 0 \\
}
\]
This induces a diagram
\[
\xymatrix{
K_0(M) \ar[r]\ar[d] & K_0(N) \ar[r] \ar[d] 
& K_0(\beta M\times S^1 \times \bR) \ar[d] \\
K_0(\pt)\ar[r] & K_0(\bC)\ar[r] &\,K_0(\bC\smallsetminus \{0\})\,\\
}
\]
with $[D_E]\in K_0(N)$ mapping by the upper right horizontal arrow
to the class of its restriction to the open subset $\beta M\times S^1
\times \bR $. Now the bottom row of this diagram is part of the exact
sequence 
\begin{multline*}
K_1(\bC)=0\to K_1(\bC\smallsetminus \{0\})\cong\bZ \\ \to
K_0(\pt)=\bZ \to K_0(\bC)\cong\bZ \to 
K_0(\bC\smallsetminus \{0\})\cong\bZ \to K_{-1}(\pt)=0.
\end{multline*}
From this we see that the map $K_0(\bC) \to K_0(\bC\smallsetminus
\{0\})$ is an isomorphism, and $K_0(\pt)\to K_0(\bC)$ is the
$0$-map. It follows that $\inda([D_E])$ can be identified with the
image of the restriction of $D_E$ to $\beta M\times S^1 \times \bR$ in
the group $K_0(\bC\smallsetminus \{0\})\cong\bZ$. But by the
assumptions on $D$, this restricted operator splits as the external
product 
\[
\left[{\left(D_{\beta M}\right)}_{E_1}\right] \boxtimes
\left[D_{S^1}\right]\boxtimes\bigl[D_\bR\bigr],
\]
and maps simply to
$k\boxtimes\left[D_{S^1}\right]\boxtimes\bigl[D_\bR\bigr]$ in
$K_0(\bC\smallsetminus \{0\})=K_0(\pt \times S^1 \times \bR)$,
where $k$ is the image of $\left[{\left(D_{\beta
M}\right)}_{E_1}\right]$ in $K_0(\pt)$. But this in turn is just 
$\ind{\left(D_{\beta M}\right)}_{E_1}$,
while $\left[D_{S^1}\right]\boxtimes\bigl[D_\bR\bigr]$ is the
canonical generator of $K_0(S^1\times \bR)$.
\end{demo}
\begin{rem}
This theorem is in many respects unsatisfactory, since it ignores what
happens on $\interior M$, and also since it fails to take advantage of
the free $S^1$-action on $\partial M$. One seemingly obvious
alternative would be to work with the $\bR$-action on $N$ which is
trivial on $M$ and which is defined on $\beta M\times S^1\times
[0,\infty)$ by the formula
\[
t\cdot(x,e^{i\theta},s) = (x,e^{i(\theta+ts)},s).
\]
This captures all the $\eta$-structure on $M$, but unfortunately it
leads to exactly the same index invariants as the ones we've already defined,
because of the fact that the forgetful map $KK^\bR\to KK$ is an
isomorphism (\cite{KaspConsp}, \S5, Theorem
2).  Another possibility would be to use the
$S^1$-action on $N'$. (There is no continuous $S^1$-action on $N$
itself.) Let $G=S^1$ and $R=R(G)=\bZ[t,\,t^{-1}]$.
We still hope to define out of $D_E$ a class in $KK^G(C_0(N'), \bC)$. 
However, existence of such a class seems to be a delicate matter since
there is no obvious reason why the restriction of $D_N$ (initially
defined on $C^\infty$ spinors on all of $N$) to those spinors whose
restriction to $\partial M= \beta M\times S^1$ is constant in the
$S^1$-factor should be essentially self-adjoint.  But if this were the
case, or if one could substitute some suitably modified operator, it
should define a $G$-equivariant Kasparov class on $N'$.  Then the 
\emph{analytic $\eta$-index} $\inda(D_E)$ of $D_E$ would again be
defined to be $c_*([D_E])$, or in other words the Kasparov product of
$[D_E]\in KK^G(C_0(N'), \bC)$ with the class of the
$G$-equivariant homomorphism
\[
c^*\co C_0(\bC)\hookrightarrow C_0(N').
\]
It would take its values in $K_0^G(\bC)\cong R=\bZ[t,\,t^{-1}]$.
(See Proposition \ref{prop:etaalgKthy}.)

Note that if $\fp$ denotes the augmentation ideal $(t-1)$ of
$R=\bZ[t,\,t^{-1}]$, then by the Localization Theorem
(\cite{Segal}; see also \cite{RosGSig}, \S3 for slight
generalizations)
$K_0^G(N')_\fp \cong
K_0^G({N'}^G)_\fp \cong K_0(M_\Sigma)\otimes R_\fp$
and $K_0^G(\bC)_\fp \cong K_0(\bC^G)_\fp \cong K_0(\pt)
\otimes R_\fp$, and $c_*$ localized at $\fp$ really can be identified
with the result of the collapse map $M_\Sigma\to \pt$.
Thus, after localizing, we would recover the original ``naive'' notion
of index theory on a $\Sigma$-manifold, where an elliptic operator
gives a class in $K_0(M_\Sigma)$ and the index is obtained via
the collapse map to a point.
\end{rem}
\begin{thm}[real $\eta$-manifold index theorem]\label{thm:realetaindthm}
Let $(M, \phi)$ be a closed spin $\eta$-manifold of
dimension $i$, and let $E$ be a real vector bundle on $M$ whose 
restriction to $\partial M$ is
pulled back from a bundle $E_1$ on $\beta M$. Fix a Riemannian
structure on $M$ compatible with its $\Sigma$-structure, and let $D_E$
be the Dirac operator on $M$ with coefficients in $E$, computed with
respect to a connection on $E$ whose restriction to a
neighborhood of the boundary is pulled back from
$\beta M$. Then the real analytic index of $D_E$ is given by
\[\inda(D_E)=
\ind\left({\left(D_{\beta M}\right)}_{E_1}\right)
\in KO_{i-2}(\pt)
\]
the index of the real Dirac operator on $\beta M$ with coefficients in 
$E_1$, which in turn is computed by applying the usual
{\lp}real{\rp} Atiyah-Singer 
Theorem {\lp}see \textup{\cite{LawM}, Chapter III, \S16}{\rp} on $\beta M$.
\end{thm}
\begin{demo}The proof is almost identical to that of Theorem
\ref{thm:etaindthm}. We only indicate the differences:

1. The only case where the map $KO_i(\pt)\to KO_i(\bC)
\cong KO_{i-2}(\pt)$ could possibly be non-zero is when $i\equiv 4
\mod 8$. In this case we have the exact sequence
\begin{multline*}
KO_4(\pt)\cong\bZ \to KO_4(\bC)\cong KO_2(\pt)\cong \bZ/2 \\
\to 
KO_4(\bC\smallsetminus \{0\})\cong KO_3(S^1) \stackrel{\partial}{\to}
KO_3(\pt)=0,
\end{multline*}
which again shows that the map $KO_4(\pt)\to KO_4(\bC)$ is zero.
So in all cases, $KO_i(\bC)$ injects into 
$KO_i(\bC\smallsetminus \{0\})\cong KO_{i-i}(S^1)\cong KO_{i-i}(\pt)
\oplus KO_{i-2}(\pt)$, in fact as a direct summand.

2. In the real case, the distinction between the two spin
structures on $S^1$ becomes relevant. Since we are using the
non-bounding spin structure, $[D_{S^1}]\in K_1(S^1)\cong \bZ$ from the proof
of Theorem \ref{thm:etaindthm} has to be replaced by $[D_\eta] \in
KO_1(S^1) \cong \bZ \oplus \bZ/2$. Note that $[D_\eta]$ projects to
the generator of $\bZ/2$ in $KO_1(\pt)$, while $[D_{S^1}]$, the Dirac
operator for the bounding spin structure, projects to $0$ in this
factor. But the two Dirac classes have the same projection in $\bZ
= \widetilde{KO}_1(S^1) \cong KO_0(\pt)$,
since they differ by the action of $H^1(S^1;\bZ/2)\cong \bZ/2$.
Thus the distinction between $[D_{S^1}]$ and $[D_\eta]$ turns out not
to matter after all, as both 
$\left[D_{S^1}\right]\boxtimes\bigl[D_\bR\bigr]$
and $\left[D_{\eta}\right]\boxtimes\bigl[D_\bR\bigr]$
project to the same thing in the image of
$KO_i(\bC)$ in $KO_i(\bC\smallsetminus \{0\})\cong KO_i(S^1\times\bR)$.
\end{demo}

\section{Applications to \psc}\label{sec:psc}
In this section we illustrate the use of the index theorems of Section
\ref{sec:index} by reproving some of the results of \cite{BotvPSC}
on \psc. First, a simple definition:
\begin{defn}\label{def:psc}
Let $(M,\phi\co \partial M\to \beta M\times P)$ be a manifold with
singularities $\Sigma=(P)$, as in Definition \ref{def:manwithsing}.
We assume $P$ is equipped with a standard scalar-flat metric.
(In our cases, $P$ will be $S^1$ or $\bZ/k$, so this will simply be
the usual metric on $P$.)
A \emph{metric of \psc} on $M$ means a Riemannian metric on $M$ which
in a collar neighborhood of the boundary diffeomorphic to
$\beta M\times P\times [0,\varepsilon)$ is a product metric
of the form
\begin{center}
(metric on $\beta M$) $\times$ (standard metric on $P$)
\\ $\times$ (standard metric on $[0,\varepsilon)$)
\end{center}
and which has {\psc} everywhere. Note that since $P\times
[0,\varepsilon)$ is scalar-flat, existence of such a metric implies
that $\beta M$ admits a metric of \psc.
\end{defn}
The following problem, first treated in \cite{BotvPSC}, now
arises: 
\begin{quest}\label{quest:psc}
Suppose $\beta M$ admits a metric of \psc. Then does $M$ admit
a  metric of {\psc} in the sense of Definition \ref{def:psc}?
\end{quest}
In this regard we have the following result:
\begin{thm}\label{thm:pscobstr}
Let $M^n$ be a closed spin $\bZ/k$-manifold. Then if $M$ admits
a  metric of {\psc} in the sense of Definition \ref{def:psc}, the
analytic index of the Dirac operator of $M$ {\lp}in the sense of
Definition \ref{def:analZkind}{\rp}
must vanish in $K_n(\pt;\bZ/k)$ or $KO_n(\pt;\bZ/k)$.
\end{thm}
\begin{demo}Assume $M$ is a $\bZ/k$-manifold with a metric of {\psc},
and form the (complex or $\cC\ell_n$-linear real) Dirac
operator $D$ with respect to this particular choice of metric. Let $N$ be
$M$ with a half-infinite cylinder attached, as in Definition
\ref{def:CaZk}. Since $N$ is complete,
the Lichnerowicz identity $D^2=\nabla^*\nabla + \frac{s}{4}$
(see \cite{LawM}, Chapter II, Theorem 8.8) is an equality of
self-adjoint operators. Here $s$ is the scalar curvature function of
$N$, and since $M$ is compact and $s$ is positive on $M$ and
translation-invariant on $\partial M\times [0,\infty)$, 
$s$ is uniformly bounded below on $N$ by a positive
constant. This implies the partial isometry part $U$ of the polar
decomposition of $D$ is unitary.  Let $\cH=\cH_0\oplus
\cH_1$ be the Hilbert space of $L^2$-spinors, on which $D$ acts. 
Note that $U$ is of the form
$\begin{pmatrix}0&U_0^*\\U_0&0\end{pmatrix}$,
where $U_0$ is a unitary operator from $\cH_0$ onto $\cH_1$. The 
Kasparov class $[D]$ is defined by $\cH$,
$U$, the action of $C^*(M;\bZ/k)$, and perhaps an additional action of
a Clifford algebra.

Now we claim that $\inda D=0$. Let $A= C^*(\pt;\bZ/k)\cong
\{f\in C_0([0,\infty), M_k) \mid f(0) \text{ a multiple of } I_k\}$,
and let $B$ be the larger unital algebra $\{f\in C_0([0,\infty]$, 
$M_k)  \mid f(0) \text{ a multiple of } I_k\}$. (Scalars here may be either
real or complex, depending on the context.) Then $B$ also acts on $\cH$
and $(B,\cH, U)$ defines a class in $K^{-n}(B)$, i.e., our class in
$K^{-n}(A)$ lies in the image of a class in $K^{-n}(B)$ defined by the same
operator $U$. The reason is that $C^\infty$ functions $\varphi$ on $N$
that are eventually constant on each cylinder $\beta M\times [0,\infty)$
have vanishing gradient in a neighborhood of infinity, and thus their 
commutator with $D$ has compact support. On the other hand, $U=D\vert
D\vert^{-1}$ and $\vert D\vert^{-1}$ is pseudodifferential of
order $-1$. From this one can deduce $[U,\varphi]=[D,\varphi]\vert D\vert^{-1}
+ D [\vert D\vert^{-1}, \varphi]$ is compact; the proof of this is
quite similar to Proposition 3.3 in \cite{HigsonZkind}.
Here are the details.  The first term, $[D,\varphi]\vert D\vert^{-1}$,
is the product of a negative-order pseudodifferential operator with
multiplication by a function of compact support, so this is
compact. The second term is pseudodifferential of negative order and
hence bounded; we want to show it is compact. We have
(using an identity from the proof of \cite{Black}, Proposition
17.11.3):
\begin{equation}\label{eq:integralexpansion}
D [\vert D\vert^{-1}, \varphi] = \frac{1}{\pi}\int_0^\infty
\lambda^{-1/2}D\left[ (D^2+\lambda)^{-1}, \varphi\right]\, d\lambda,
\end{equation}
but
\[
[ (D^2+\lambda)^{-1}, \varphi] = - (D^2+\lambda)^{-1}[(D^2+\lambda),
\varphi](D^2+\lambda)^{-1}
\]
and $[(D^2+\lambda),\varphi]=[D,\varphi]D + D[D,\varphi]$ has compact
support, hence 
\[
[(D^2+\lambda), \varphi](D^2+\lambda)^{-1}
\]
is compact, and then
\[
D\left[ (D^2+\lambda)^{-1}, \varphi\right]=
-D  (D^2+\lambda)^{-1} \left[ (D^2+\lambda),
\varphi \right] (D^2+\lambda)^{-1}
\]
is compact. Then by equation (\ref{eq:integralexpansion}),
the remaining term in $[U,\varphi]$ is also compact.
Since $U^2=1$ and $U=U^*$, there is nothing to
check as far as the other axioms for a Kasparov bimodule are 
concerned, so $(B,\cH, U)$ defines a class in $K^{-n}(B)$.
Furthermore, we have a short exact sequence
\[
0\to C_0((0,\infty], M_k)\to B\to C(\pt)\to 0,
\]
with the ideal $C_0((0,\infty], M_k)$ contractible, and thus
$K^{-n}(B)\cong K_n(\pt)$. 
The map $B\to C(\pt)$ is split by the inclusion of scalar multiples
of the identity, under which the class of $(B, \cH, U)$ pulls back
to the class of $(\cH, U)$.  So via the isomorphism
$K^{-n}(B)\cong K_n(\pt)$, we see that the
class of $(B,\cH, U)$  can be identified with class of
$U$ (perhaps with some auxiliary Clifford algebra action, if we're in
the case of the $\cC\ell_n$-linear real Dirac
operator), which vanishes, since $U$ is unitary. Hence $c_*([D])=0$.
\end{demo}
This gives another proof of the ``only if'' direction of the following
theorem from \cite{BotvPSC}:
\begin{thm}[Botvinnik \cite{BotvPSC}]\label{thm:botvZtwo}
Let $M^n\!$ be a closed spin $\bZ/2$-manifold, with $n\ge 6$,
and assume $M$ and $\beta M$ are connected and simply connected.
Then $M$ admits a metric of {\psc} in the sense of Definition 
\ref{def:psc} if and only if the image $\alpha_\Sigma(M)$ of the
canonical class defined by the spin structure in $KO_n(M;\bZ/2)$
vanishes in $KO_n(\pt;\bZ/2)$.
\end{thm}
\begin{rem}Note that the obstruction $\alpha_\Sigma(M)\in KO_n(M;\bZ/2)$
of Theorem \ref{thm:botvZtwo} includes within it the obstruction to
existence of a {\psc} metric on $\beta M$. Indeed, since $n\ge 6$,
$\dim \beta M\ge 5$, so we know from \cite{Stolz} that $\beta M$
admits a metric of {\psc} if and only if the usual index invariant
$\alpha(\beta M)\in KO_{n-1}(\text{\textup{pt}})$ vanishes. Now
the existence of spin structure on $M$ with $\partial M\equiv
\beta M\amalg \beta M$ implies that the class of
$\beta M$ must be a $2$-torsion class in $\Oms_{n-1}$, which forces
the $\widehat A$-genus of $\beta M$ to vanish if $n\equiv1$ (mod $4$).
On the other hand, we claim that under the exact sequence
\begin{multline*}
\cdots \longrightarrow KO_n(\text{\textup{pt}}) \stackrel{2}{\longrightarrow}
KO_n(\pt) {\longrightarrow} KO_n(\pt;
\bZ/2) \\ \stackrel{\partial}{\longrightarrow}
KO_{n-1}(\text{\textup{pt}}) \stackrel{2}{\longrightarrow} \cdots,
\end{multline*}
$\partial(\alpha_\Sigma(M))=\alpha(\beta M)$, so $\alpha_\Sigma(M)\ne 0$
if $\alpha(\beta M)$ is a non-zero $2$-torsion class. 
For example, if $\beta M$ is an exotic $9$-sphere with
$\alpha(\beta M)$ a non-zero $2$-torsion class in $\Oms_9$, then
$\beta M\amalg \beta M$ bounds a spin $10$-manifold $M$ which can be
given a $\bZ/2$-manifold structure, and $\alpha_\Sigma(M)\ne 0$.

To check this, observe that by Theorem \ref{thm:Zkindthm}, the
topological invariant $\alpha_\Sigma(M)$ coincides with $\inda(D)$,
$D$ the Dirac operator on $M$ (or more exactly on $N$, the manifold
with cylinders attached). Let $[D]\in KO^{-i}(C^*(M;\bZ/2))$ 
be the associated class. Via the exact sequence in Proposition
\ref{prop:Zkalgstructure}, this restricts to $[D_{\beta M}]$, the class
of the Dirac operator on $\beta M$, in 
\[
KO^{-i}\bigl(C_0(\bR)\otimes C(\beta M) \otimes M_2\bigr)
\cong KO_{i-1}(\beta M). 
\]
\end{rem}
Similarly, Theorem \ref{thm:pscobstr} is related to the
``only if'' direction of the following theorem from \cite{BotvPSC}:
\begin{thm}[Botvinnik \cite{BotvPSC}]\label{thm:botveta}
Let $M^n$ be a closed spin $\eta$-manifold, with $n\ge 7$,
and assume $M$ and $\beta M$ are connected and simply connected.
Then $M$ admits a metric of {\psc} in the sense of Definition 
\ref{def:psc} if and only if the image $\alpha_\Sigma(M)$ of the
canonical class defined by the spin structure in $KU_n(M)$
vanishes in $KU_n(\pt)$.
\end{thm}
\begin{thm}\label{thm:psceta}
Let $(M^n, \phi)$ be a closed spin $\eta$-manifold, and suppose
$M$ admits a metric of {\psc} in the sense of Definition 
\ref{def:psc}. Then the analytic $\eta$-index
$\inda(D)\in KO_n(\bC)\cong
KO_{n-2}(\pt)$ of Dirac operator on $M$ {\lp}as defined
in Definition \ref{def:analetaind}{\rp} must vanish.
\end{thm}
\begin{demo}By the index theorem \ref{thm:realetaindthm},
$\inda(D)$ may be identified with the index of the real Dirac
operator on $\beta M$, which is an obstruction to positive
scalar curvature on $\beta M$ \cite{Hit}. However, a metric
of {\psc} on $M$ must by definition restrict to a metric of {\psc}
on $\partial M\cong \beta M\times S^1$ which is a product of a metric
on $\beta M$ with the standard flat metric on $S^1$. So such a metric
can only exist when $\beta M$ admits a  metric of {\psc}.
\end{demo}
\begin{rem}One can also give another proof of Theorem \ref{thm:psceta}
along the lines of the proof of Theorem
\ref{thm:pscobstr}.  In addition, the argument used to prove Theorem
\ref{thm:pscobstr} can be extended to give an interpretation of the
$\bZ/k$-index of the Dirac operator in a more general context.
\end{rem}
\begin{thm}\label{thm:ZkindofDirac}
Let $M^n$ be a closed spin $\bZ/k$-manifold, equipped with a
$\bZ/k$-metric restricting on $\beta M$ to a metric of \psc. Let
$N=M\cup_{\partial M}\partial M\times [0, \infty)$ and let $D$ be the
$\cC\ell_n$-linear real Dirac operator on $N$. Then $D$ {\lp}acting on the
$\bZ/2$-graded Hilbert space $\cH$ of $L^2$ sections of the
appropriate bundle of free right $\cC\ell_n$-modules
over $N${\rp} has finite-dimensional kernel, and
the $\bZ/k$-index of $D$ is the mod-$k$ reduction of the
$KO(\pt)_n$-valued index of $D$ {\lp}computed from the kernel of $D$
over $N$, viewed as a graded $\cC\ell_n$-module{\rp}.
\end{thm}
\begin{demo}Because $N$ has \emph{uniformly} {\psc} on the ends $\beta
M\times [0, \infty)$, the spectrum of $D$ is bounded away from $0$ on
the complement of a finite-dimensional subspace of $\cH$
(\cite{GromovLawson}, Theorem 3.2). Then if $U$ is the partial
isometry part of the polar decomposition of $D$, $U$ is Fredholm and
the $\bZ/k$-index of $D$ is defined by the triple $(C^*_\bR(\pt; \bZ/k),
\cH, U)$. As in the proof of \ref{thm:pscobstr}, this Kasparov class
is in the image of $KO^{-n}(B)$, where $B$ is a suitable unital
extension of $C^*_\bR(\pt; \bZ/k)$ with the same $KO$-theory as the
scalars.  The class in $KO^{-n}(B)$ that maps to the $\bZ/k$-index of
$D$ is defined by $(U, \cH)$ (together with the $\cC\ell_n$-action),
and the map 
\[
KO^{-n}(B)\to KO^{-n}(C^*_\bR(\pt; \bZ/k))
\]
can be identified with the map $KO_n(\pt)\to KO_n(\pt; \bZ/k)$, which is
simply reduction mod $k$. So $\inda D$ is the mod $k$ reduction of the 
$KO(\pt)_n$-valued index of $D$.
\end{demo}

\bibliography{sing}
\bibliographystyle{amsplain}

\vspace*{15pt}
\begin{flushleft}
Department of Mathematics\\
University of Maryland\\
College Park, MD 20742-4015\\
email: \texttt{jmr@math.umd.edu}
\end{flushleft}
\end{document}